\documentclass{article}
\usepackage[affil-it]{authblk}
\usepackage{amssymb,amsmath}
\usepackage[pdftex]{graphicx}
\usepackage{hyperref}

\usepackage[utf8]{inputenc}
\usepackage[english]{babel}

\usepackage{natbib}
\makeatletter
\def\fps@figure{htbp}
\makeatother

\date{}

\begin{document}

\title{A Simple 3D Isometric Embedding of the Flat Square Torus}

%\author{
%J. Richard Gott III\footnote{Department of Astrophysical Sciences,
%Princeton University,
%Princeton, NJ, 08544, USA } 
%\and
%Robert J.  Vanderbei\footnote{Department of Operations Research and Financial
%        Engineering,
%Princeton University,
%Princeton, NJ, 08544, USA }
%}

\author[1]{
J. Richard Gott III\thanks{Corresponding author: jrg@astro.princeton.edu}
}
\author[2]{
Robert J.  Vanderbei
}

\affil[1]{Department of Astrophysical Sciences %,
%Princeton University,
%Princeton, NJ, 08544, USA 
} 

\affil[2]{Department of Operations Research and Financial
        Engineering %,
Princeton University,
Princeton, NJ, 08544, USA }

\maketitle

\begin{abstract}
Start with \cite{Gott}'s envelope polyhedron (Squares---4 around a
point): a unit cube missing its top and bottom faces. Stretch by a
factor of 2 in the vertical direction so its sides become ($2\times1$ unit)
rectangles. This has 8 faces (4 exterior, 4 interior), 8 vertices, and
16 edges. $F - E + V = 0$, implying a (toroidal) genus = 1. It is
isometric to a flat square torus. Like any polyhedron it has zero
intrinsic Gaussian curvature on its faces and edges. Since 4 right
angled rectangles meet at each vertex, there is no angle deficit and
zero Gaussian curvature there as well. All meridian and latitudinal
circumferences are equal (4 units long).
\end{abstract}

\section{Introduction}

Many video games including \emph{Pac-Man} and \emph{Asteroids} are
played on a square TV screen with a toroidal geometry. Pac-Man, or a
spaceship that disappears off the top of the screen immediately
reappears at the bottom at same horizontal location. Likewise, if
Pac-Man or a spaceship disappears off the right-hand edge of the screen
it immediately reappears on the left-hand edge of the screen at the same
vertical location. These games thus have a square Euclidean flat
geometry but a toroidal topology. It takes seven colors in general to color nations
on such a TV game map, rather than the four required on the plane. It is
a map of a toroidal flatland universe of two dimensions, that obeys
Euclidean geometry, has zero Gaussian curvature everywhere (every small
circle around any point has a circumference of $2 \pi r$, and the sum of
angles in any triangle on the screen has a sum of angles of $180^{\circ}$). But
the map has two boundary cuts, one across the top of the screen, one
across the right of the screen, connected to the same cuts on the
opposite sides. This game can be mapped onto a donut (or tire inner tube
shape) in 3D space, but with distortion, since the tire has a curved
surface. It has positive Gaussian curvature on the outer circumference
of the tire, and negative Gaussian curvature on the interior of the
tire, where the geometry is saddle shaped.{~ }Can the game geometry be
embedded in 3D space? Videos illustrating this usually start with a
square, flat sheet of paper. The demonstrator then bends it into a
cylinder and tapes it together. If the side length of the square is $S$,
the cylinder has a height of $S$ and a radius of $r$, where $S = 2 \pi r$. The
demonstrator bends the paper; this does not change its intrinsic
Gaussian curvature of zero everywhere. A cylinder has extrinsic
curvature, but zero intrinsic curvature everywhere. One of the two
boundary cuts has been sown together and healed when the demonstrator
tapes the two opposite edges together. This leaves at the ends of the
cylinder, two circular edges, which now need to be taped together. Then
the demonstrator tries fruitlessly to bend the cylinder around to bring
them together, and illustrates that it would tear and crumple the
paper to try to do so. One is trying to produce a tire inner tube with
curvature from a flat cylinder with no curvature. Doesn't work we are
told.

But if we lived in 4D Euclidean space with a metric
$ds^2 = dx^2 + dy^2 + dz^2 + dw^2$, we could do it by bending
the cylinder around smoothly in the $w$ direction. The result is a 2D
flat torus living in 4D space (a Clifford torus):
\[ 
        r^2 = x^2 + y^2 = \frac{1}{2} = z^2 + w^2 .
\]
This is a 2D surface defined by two equations, living in 4D Euclidean
space.{~}

Here $r$ is a constant $r = 1/\sqrt{2}$. Define $\theta$ and $\phi$ by 
$x = r \sin({\theta})$, $y = r \cos({\theta})$, 
$z = r \sin( \phi )$, $w = r \cos( \phi )$ and we can show that the
equations of the surface are satisfied automatically since
$\sin^2 + \cos^2 = 1$, and that using $(\theta, \phi)$ 
as coordinates on the surface, the metric on the surface is
\[
ds^2 = r^2 d\theta^2 + r^2 d\phi^2, \qquad 
\text{where $\quad 0 \le \theta < 2 \pi, \;\;\; 0 \le \phi < 2 \pi$}.
\]
This is a flat metric with zero Gaussian curvature. It is a flat square
with side $S = 2 \pi r = \sqrt{2} \pi$, our original TV game with no distortion, and no
boundary cuts. The vertical and horizontal coordinates are cyclic, and
the boundary cuts have been healed. This is an isometric embedding of
the TV game geometry in 4D Euclidean space. 
%It is called an isometric
%embedding because all metric distances between points in the TV game and
%their image points in the toroidal surface in 4D Euclidean space are
%preserved.

\cite{Nash} famously proved that there must exist a
${\mathbb C}^{1}$ 
isometric embedding of the flat square torus
in 3D Euclidean space -- but he provided no constructive proscription
for doing so. (A ${\mathbb C}^{1}$ isometric embedding has
continuous first derivatives, in other words it admits tangent planes
that very smoothly, being defined at every point. (A
${\mathbb C}^{1}$ isometric map which is not
${\mathbb C}^2$ has no defined extrinsic curvature
everywhere.) Recently, Borrelli, Jabrane, Lazarus, and Thibert in
their paper ``Isometric Embedding of the Square Flat Torus in Ambient
Space'', \cite{Borrelli} have shown how to do it. It looks like a torus, but with ripples
on it. And there are ripples on those ripples tipped at an angle, and
tinier ripples on those ripples at another angle, ad infinitum. It is
sort of a smooth fractal.{~ }The amplitude of the ripples on the inner
edge of the torus is larger, making its circumference exactly as large
as the smoother circumference on the outside. Beautiful pictures of it
are shown in their paper. For his original work in this field, John Nash
(who by the way often visited our Astrophysics Library in Princeton)
received the Abel Prize in Mathematics. Many mathematicians regard this
as a greater work on his part than his work on game theory for which he
won the Nobel Prize.{~}

While the embedding of Borelli, et al. is beautiful and important
mathematically, it does not make a very practical object on which to
play the game. It is so crinkled that it would be hard to even see the
game portrayed on its corrugated surface. Yet it is isometric.

\section{Envelope Polyhedra}

\begin{figure} 
\includegraphics[width=5in]{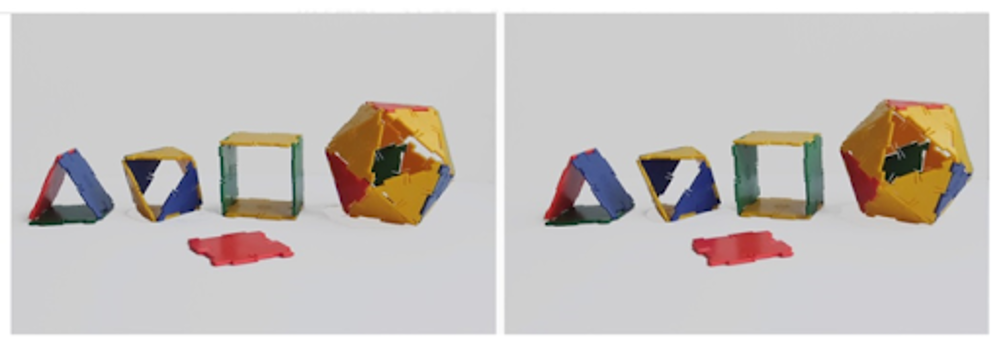}
\caption{
Some Finite Envelope Polyhedra from \cite{Gott}. Envelope
polyhedra are regular polyhedra with regular polygons as faces, the
arrangement of polygons around each vertex must be identical, but not
all dihedral angles are equal, and some are $0^{\circ}$, allowing polygons to
appear back to back. From left to right at top: Squares---4 around a
point, Triangles---6 around a point, Squares---4 around a point (a unit
cube missing its top and bottom---tipped on its side in this picture---a
subject of this paper), and Triangles---8 around a point. At bottom,
Squares---2 around a point (a dihedron--see \cite{Coxeter}). To
see this stereo view, touch your nose to the page and slowly withdraw
it. You will see a blurry 3D image in the center which will come into
focus.{~}
}
\label{fig1}
\end{figure}

We wish to propose here a simple ${\mathbb C}^{0}$
isometric mapping of the flat square torus into 3D Euclidean space. Note
that in Figure \ref{fig1}, Gott's envelope polyhedron Squares---4 around a point
has a toroidal topology. It is a cube with its top and bottom
missing.{~}

The cube with its top and bottom faces missing has 4 exterior squares linked by the square holes at the top and
bottom to 4 interior squares. It thus has 8 faces, 4 outside edges, 4
inside edges, and 8 edges at the top and bottom that link the exterior
and interior faces. It has 8 vertices. 
Thus $F - E + V = 8 - 16 + 8 = 0 = 2(1-g)$, giving this a genus of 1, 
which is a torus. Notice that there is
no Gaussian curvature on any of the faces, they are flat. There is no
Gaussian curvature on any of the edges; they are limits of flat surfaces
of partial cylinders whose radii shrink to zero. As in any polyhedron,
the Gaussian curvature is contained in delta-functions at the vertices.
In this case since there are 4 squares around every vertex, the angle
deficit is zero at each vertex and the integral of the Gaussian
curvature at the vertex is zero. This is a geometric surface which
contains zero volume. The circumference around its equator on the
outside is 4 units (if the cube is of one-unit side length). This
circumference circles the 4 outside faces. It is obviously a geodesic.
The inner circumference circling the 4 inside faces is also 4 units. It
is also a geodesic. What about the upper edge: the four straight edges
at the top that form the opening that connects the outer and inner faces
at the missing top face of the cube. These seem to turn $90^{\circ}$ four times
as one circles the top of the cube. Yet this circuit is also a geodesic.
For any 2D surface imbedded in 3D Euclidean space, the shortest curve
connecting points $A$ and $B$ in the surface is a geodesic. Pick any two
points $A$ and $B$ in the four top edges. The shortest path from $A$ to $B$
staying in the surface will necessarily travel along those top edges. If
the curve dips down into any of the faces during part of its journey it
will be longer.

Take this cube and lengthen it in the vertical direction by a factor
of 2. Now there is a missing square at the top and bottom, and there are
4 rectangular sides made up of 2:1 rectangles. There are 4 exterior
rectangles of length 2 units and width 1 unit, and 4 interior rectangles
of length 2 units and width 1 unit. The equatorial circumference of any
horizontal geodesic traversing the 4 exterior rectangular faces is 4
units as before. But now consider the longitudinal circumference of the
torus. Starting at a bottom edge, one first traverses an outside
rectangular face traveling 2 units, one then travels over the top edge
and traverses the inside rectangular face traveling another 2 units and
reaches the bottom edge where one started. These are now the vertical
lines in the original game, 4 units long, while the equatorial
geodesics, also 4 units in circumference, are the horizontal lines in
the original game. This is a flat square torus.{~}

\section{Origami}

Go back to the original demonstration. Bend the original square
piece of paper containing an \emph{Asteroids} game into a cylinder with
the game on the outside.{~ }Tape it together to form a complete cylinder
as before: the top of the original square piece of paper is now taped to
the bottom. Now instead of giving up, we are solving the problem with
origami.

Origami has an interesting history in mathematics. In 1980, Hisashi
Abe solved the angle trisection problem by using origami: the angle
trisection takes 7 steps, see for example \cite{Alperin} and
\cite{Fuchs} for discussion, proofs and the axiomatic basis for
this.{~ }It didn't occur to the ancient Greeks to think of a plane as a
piece of paper that could be folded. Likewise, the ${\mathbb C}^{1}$
embeddings of \cite{Nash} and \cite{Borrelli}
by definition do not allow folds, although it is only
the intrinsic geometry of the embedded surface that we are really
interested in.

Now take the horizontal cylinder you have just formed lying on a
table, where the taping line is flush with the table, and squash the
cylinder flat: crease it making two folds, opposite each other on the
cylinder and equidistant from the taping line. Let the square have side
length of 4 units. One now has a flattened, open-ended envelope (4 units
long and 2 units wide) with the game printed on the top and bottom
outside sides of this envelope.  Flip it over so the taped edges are on the top
side. Now pick up the two open ended ends
(which are 2 units wide) and 4 units apart and fold them inward:{~ }fold
it in quarters so that it makes 4 rectangles, each 2 units by 1 unit.
Mate the two open-ended slits and tape them together to form a tall,
square, open-ended box with 4 rectangular sides, each 2 units tall and 1
unit wide. This is the elongated box-like structure we have just
discussed above. See Figures \ref{fig2} and \ref{fig3} below. We have made the
open-ended boxlike surface out of
our original flat square and we have taped the opposite sides of the
square together as required to make the square torus. The spaceship
shoots an asteroid with a ray gun.

\begin{figure} 
\includegraphics[width=5in]{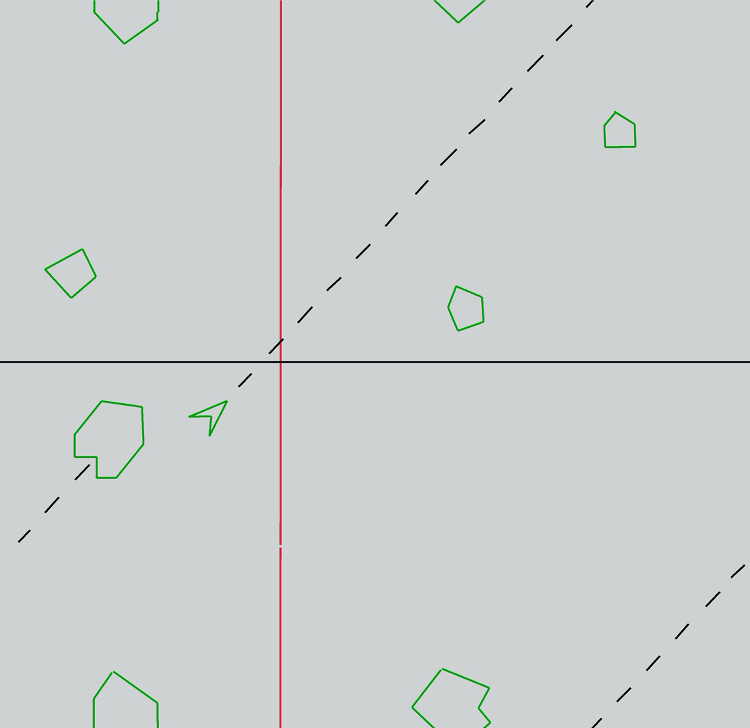}
\caption{
The square flat torus showing an \emph{Asteroids} game, top
and bottom edges are identified as are left and right edges.
}
\label{fig2}
\end{figure}

\begin{figure} 
\begin{center}
\includegraphics[width=2.853in]{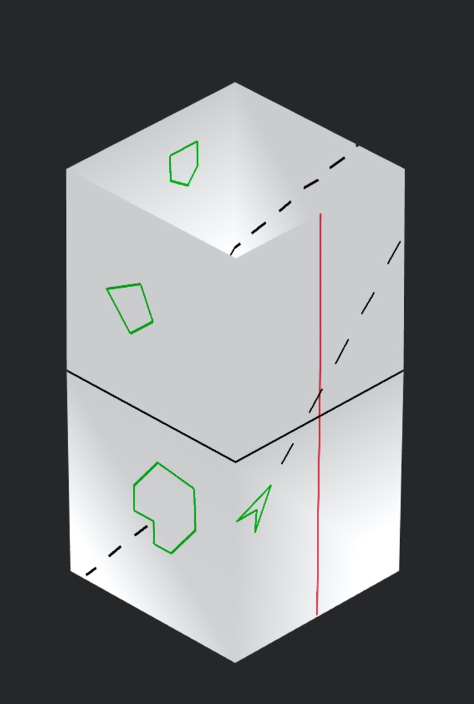}
\end{center}
\caption{
Folding the square flat torus into a simple boxlike surface
isometrically embedded in 3D space that is ${\mathbb C}^{0}$.
}
\label{fig3}
\end{figure}

A literature search revealed that people have come close, but
failed, to find this solution earlier. The \cite{Zalgaller} paper,
``Some Bendings of a Long Cylinder'' came up with some polyhedron
solutions for a long rectangular torus. These enclose volume. Basically,
he makes a long triangular prism by origami first, makes origami fold
crimps in it and ends up with an n-gonal toroidal polyhedron, where the
sum of face angles around each vertex is $360^{\circ}$. As part of his
introduction, he gives a construction where he takes a long rectangle
and folds in in half and tapes it together to make a long, ironed
sleeve. Then he takes the open cuff ends and creates two folds in the
sleeve so that he can bring the two cuff ends back together at the
center and tape them together in the center. Now he has a four-folded
surface ironed flat in the plane. This is really a long rectangular flat
torus. It has a disadvantage that half of its map of the video game is
buried on inside faces of the ironed flat surface of the plane so that
you can't see them. He calls this a direct flat torus. He seems to not
realize that this construction could be used to make a flat square torus
embedded in a 2D Euclidean space as well, instead of the long one he is
making.{~}

There is a wonderful YouTube video uploaded on June 8, 2015
by \cite{Segerman} titled ``Hinged Flat Torus''
showing him folding a flat polyhedron to form a
torus. It is a bending of a long cylinder as well. It has 10 rectangles
around its outer equatorial surface forming a decagonal prism with open
top and bottom. Then a zig-zag configuration of 20 long triangles that
connect top vertices with 
%a vertex rotated by a couple of places $(72^{\circ})$
rotated vertices
on the decagonal edge below. This encloses a volume. This is
called a piecewise linear flat torus. Unfolded, it shows a row of 10
rectangles along the bottom and 20 long, tipped triangles making a
parallelogram on the top. This figure is flat and would tile the plane.
There is no Gaussian curvature at any vertex; the sum of angles around
each vertex is $360^{\circ}$ with three triangles and two rectangles around each
point. This is longer than it is wide and so does not solve the flat
torus embedding. Guy Valette saw a torus like this at Oberwolfach over 30 years ago.
Thus, we don't know who invented it.{~}

Had they allowed polygons to appear back to back with dihedral
angles of $0^{\circ}$ as in Gott's envelope polyhedra they could have used the 10
rectangles alone to have formed an open-ended decagonal prism with
rectangles meeting 4 around a point (at the top and bottom decagonal
edges). Make the rectangles tall enough (5:1) and one would have an
envelope polyhedron that was an embedding of the square flat torus.

We propose another origami construction which embeds the flat torus
in 3D space.{~ }Start with a flat square 4 units on a side. Fold it in
half, and tape the two sides together. Now bend this 4-unit long, 2-unit
wide rectangular envelope of zero thickness into a cylinder by taping
the two open ends together. One now has an open-ended cylinder 2 units
tall with a circumference of $2 \pi r = 4$ units. The equatorial circumference
is 4 units. The longitudinal circumference is also 4 units, since
starting at the bottom circular edge traveling upward on the outside of
the cylinder 2 units, crossing over the upper edge and traveling
downward 2 units one reaches the bottom edge where one started. All the
longitudinal and equatorial circumferences are geodesics. The circular
edges at the top and bottom are also geodesics, being straight lines in
the original flat square. Take the original square, tile the plane with
images of it. Pick any two points $A$ and $B$ in the original square and
connect them with a straight line connecting their closest pair of
images in the tiling and this will be the geodesic connecting them on
the surface. Such a geodesic encountering the top or bottom edge will
make an equal angle with the edge on the other side. Unfolded, the line
will be straight. This is a surface that is the limit of an $N$-gonal
prism without its top and bottom $N$-gons as $N$ goes to infinity. 
The prisms have $N$
external rectangular faces, $N$ internal rectangular faces, $N$ external
edges, $N$ internal edges, $2N$ top and bottom edges linking the internal
and external edges, and $2N$ vertices. Thus, 
%it has 
they have
$F - E + V = 2N - 4N + 2N = 0 = 2(1-g)$, so $g = 1$ 
%and it is 
and each is
a torus.

An interesting side note: starting with the same flat square torus 4
units by 4 units, one can make this out of 16 squares. These meet 4
around a point, like in a checkerboard. This can be embedded in 3D space
as we have described above: 8 squares on the inside of a 2x1x1 hollow
box, and 8 squares on the outside. This is a cube stretched vertically
by a factor of 2 missing its top and bottom. Dihedral angles at the top
and bottom edges are $0^{\circ}$. But it is also possible to embed this in 4D
Euclidean space as 16 square faces which are a subset of the square
faces of a tesseract (or hyper-cube). The tesseract has a surface
composed of 8 cubic cells. Unfold the tesseract in a cross like shape as
shown in the Salvador Dali 1954 painting Corpus Hypercubus. Four cubes
form the vertical post of this cross. Consider the four sides of each of
these four cubes, they form a tube with cylindrical topology whose
surface consists of 16 squares. When this polyhedral net is folded back
up into the tesseract, the top of the cross will meet the bottom, to
make this a torus with 16 square faces. This is a two-sided surface. It
divides the tesseract surface into two equivalent parts, like the
Clifford Torus divides the 3-sphere into two equivalent parts. On this
torus, the square faces meet 4 around a point.

\bibliography{refs}

%REFERENCES{ }
%
%Alperin, Roger C., Trisections and Totally Real Origami, arXiv:
%math/0408159v1 (2004){~}
%
%Borrelli, Vincent; Jabrane, Said; Lazarus, Francis; Thibert, Boris;
%``Isometric Embedding of the{~}
%
%Square Flat Torus in Ambient Space, Ensaios Mathematicos, 2013, Vol. 24,
%1-91.{~}
%
%Coxeter, H. S. M., ``Regular Skew Polyhedra in Three and Four
%Dimensions, and Their{~ ~ ~ ~ ~}
%
%Topological Analogues.'' Proc., London Math. Soc. 2, 1937.
%
%Fuchs, Clemens, Elemente der Mathematik, 66 (2011) 121-131.{~}
%
%Gott, J.R., Envelope Polyhedra, arXiv:1908.05395v1 {[}math.MG{]} 15 Aug
%2019
%
%Nash, J. F., ${\mathbb C}^{1}$-isometric embeddings, Ann. of
%Math, 60 (1954), no. 3, 383-396.
%
%Henry Segerman, Hinged Flat Torus, YouTube video uploaded on June 8,
%2015
%
%Zalgaller, V. A., Some Bendings of a Long Cylinder, journal of
%Mathematical Sciences, Vol. 100, no. 3,{~}
%2228-2238
%
%{{~}}{~ ~ ~}

\end{document}